\documentclass[a4paper,10pt]{article}
\usepackage[utf8]{inputenc}
\usepackage{color}

\usepackage[all]{xy}
\usepackage{amsmath,amssymb,amsfonts,amsthm,graphicx}
\usepackage[margin=1.2in]{geometry}

\newtheorem{thm}{Theorem}[section]

\newtheorem{cor}[thm]{Corollary}
\newtheorem{prop}[thm]{Proposition}
\newtheorem{df}[thm]{Definition}
\newtheorem{lm}[thm]{Lemma}

\title{A criteria for a finite permutation group to be transitive}
\author{Julian Brough}

\begin{document}
\date{}
\maketitle

\begin{center}
\small
\textit{Faculty of Mathematics, Centre for Mathematical Sciences,}

\textit{Wilberforce Road, Cambridge, England, CB3 0WA}
\end{center}

\paragraph{}
  \textit{Keywords:}

\textit{Finite transitive permutation Groups}

\normalsize
\begin{abstract}
Let $G$ be a finite permutation group on a finite set $\Omega$.
The notion of $G$ being quasi-transitive on $\Omega$ was defined by Alan Camina \cite{Camina}; in that paper conditions were established that ensured a quasi-transitive group on a finite set $\Omega$ was transitive on $\Omega$.
The aim of this paper is to validate the conjecture made in \cite{Camina}: given any group $G$, if $G$ is quasi-transitive on a finite set $\Omega$ then $G$ is transitive on $\Omega$.
\end{abstract}

\section{Introduction}

Let $G$ be a finite 2-transitive group acting on a finite set $\Omega$; a simple observation yields that for $G_{\alpha\beta}$ the pointwise stabiliser of $\alpha,\beta\in\Omega$, then $|G_{\alpha\beta}|$ constant for all $\{\alpha,\beta\}\subset \Omega$.
Thus it is natural to ask, what does it mean for the action of $G$ on $\Omega$ if $|G_{\alpha,\beta}|$ is constant for all $\alpha,\beta\in\Omega$?
Consider the case that $|G_{\alpha\beta}|=1$ for all $\{\alpha,\beta\}\subset \Omega$; no non-trivial element in $G$ fixes more than one point of $\Omega$.
Thus $G$ must act Frobeniusly or regularly on each orbit of $\Omega$; however a finite group $G$ can act Frobeniusly in at most one way and so $G$ acts Frobeniusly on at most one orbit.
Hence $|G_{\alpha\beta}|=1$ for all $\{\alpha,\beta\}\subset \Omega$, if and only if $G$ acts Frobeniusly on one orbit and regularly on the remaining orbits.
From this we see that constant $|G_{\alpha\beta}|$ does not mean $G$ has to be transitive; however if we assert that $|G_{\alpha\beta}|=t$ for all $\{\alpha,\beta\}\subset \Omega$ and $t>1$, can we now conclude $G$ has to act transitively?
To study this hypothesis Camina made the following definition.

\begin{df}\cite[Definition]{Camina}
 A finite permutation group $G$ acting on a finite set $\Omega$ is called quasi-transitive if there is a natural number $t>1$ such that $|G_{\alpha\beta}|=t$ for all two element subsets $\{\alpha,\beta\}\subset \Omega$.
\end{df}

In \cite{Camina}, Camina established conditions for a group $G$ that ensured quasi-transitive implies transitive.

\begin{prop}\cite[Proposition 6]{Camina}\label{Camina1}
 If $G$ is quasi-transitive and has a non-trivial normal soluble subgroup then $G$ is transitive.
\end{prop}

\begin{thm}\cite[Theorem]{Camina}\label{Camina2}
 If $G$ is quasi-transitive and $G_{\alpha\beta}$ is abelian for all two-element subsets $\{\alpha,\beta\}\subset \Omega$, then $G$ is transitive.
\end{thm}

In this note we validate the conjecture made in \cite{Camina}.

\begin{thm}[\textbf{Main Theorem}]\label{MainTheorem}
 Let $G$ be a finite permutation group $G$ acting quasi-transitively on a finite set $\Omega$, then $G$ acts transitively on $\Omega$.
\end{thm}

Our proof will make use of Proposition~\ref{Camina1}, Proposition~\ref{Camina2} and a classification result by Bamberg, Giudici, Liebeck, Praeger and Saxl \cite{LPS}.

\begin{thm}\cite[Theorem 1.1]{LPS}\label{LPS1}
 Every finite primitive $\frac{3}{2}$-transitive group is either affine or almost simple.
\end{thm}
\begin{thm}\cite[Theorem 1.2]{LPS}\label{LPS2}
 Let $G$ be a finite almost simple $\frac{3}{2}$-transitive group of degree $n$ on a set $\Omega$. Then one of the following holds:
 \begin{enumerate}
  \item $G$ is $2$-transitive on $\Omega$.
  \item $n=21$ and $G$ is $A_7$ or $S_7$ acting on the set of pairs of elements of $\{1,\dots,7\}$; the size of the non-trivial subdegrees is 10.
  \item $n=\frac{1}{2}q(q-1)$, where $q=2^f\geq 8$, and either $G=PSL_2(q)$ or $G=P\Gamma L_2(q)$ with $f$ a prime; the size of the non-trivial subdegrees is $q+1$ or $f(q+1)$, respectively.
 \end{enumerate}

\end{thm}

In particular, if we have a quasi-transitive group on a finite set $\Omega$, and $\Delta$ is an orbit of $G$ on which $G$ acts as a finite almost simple $\frac{3}{2}$-transitive group, then $t$ is determined by Theorem~\ref{LPS2}.

\section{The reduction}

Before establishing the main theorem, we reduce the problem to the situation in Theorem~\ref{LPS2}.
First we introduce some notation and recall some properties about quasi-transitive groups established in \cite{Camina}.

Assume a finite group $G$ acts quasi-transitively but not transitively on a finite set $\Omega$, and $|G_{\alpha\beta}|=t$ for all two element subsets $\{\alpha,\beta\}$ of $\Omega$.
As $G$ is not acting transitively there are at least two orbits of $G$; label the orbits for this action by $\Delta_1,\Delta_2,\dots,\Delta_r$ for some $r>1$.
For each $i=1,\dots,r$ fix an element $\alpha_i\in\Delta_i$; the subdegrees of $G$ on $\Delta_i$ are given by $d_i=|G_{\alpha_i}|/t$.
As $d_i$ is constant for the whole orbit $\Delta_i$; the group $G$ acts $\frac{3}{2}$-transitively on each $\Delta_i$.
In addition, by the Orbit-Stabiliser Theorem, $|G_{\alpha_i}|=\frac{|G|}{|\Delta_i|}$, thus $td_i=|G_{\alpha_i}|=\frac{|G|}{|\Delta_i|}$ and hence $d_1\cdot|\Delta_1|=d_i\cdot|\Delta_i|$ for all $i=1,\dots r$.
In \cite{Camina} Camina established arithmetical conditions between the $d_i$ and a condition for the action of $G$ on $\Delta_i$.

\begin{lm}\cite[Lemma 1]{Camina}\label{dCoprime}
 Given $d_i$ as above, the $d_i$'s are pairwise coprime.
\end{lm}

\begin{prop}\cite[Proposition 4]{Camina}\label{qtfaith}
 Let $G$ be a quasi-transitive group then $G$ acts faithfully on each orbit.
\end{prop}

In particular none of the orbits $\Delta_i$ are a singleton element by Proposition~\ref{qtfaith}.
These two results have the following Corollary.
\begin{cor}\label{distinct}
 Given $\Delta_i$ and $d_i$ as above, $d_1\cdot|\Delta_1|=d_i\cdot|\Delta_i|$ for all $i=1,\dots r$.
 In addition, the set of $|\Delta_i|$ are pairwise distinct.
 \begin{proof}
  The first result is discussed in the comments above. Hence it only remains to prove the second point.

  Assume $|\Delta_i|=|\Delta_j|$, by the first part $d_i=d_j$ and hence $d_i=1$ by Lemma~\ref{dCoprime}.
  As $t=|G_{\alpha\beta}|$ for all two element subsets $\{\alpha,\beta\}$ of $\Omega$, then $|G_{\alpha_i}|=|G_{\alpha_i\beta}|$ for all $\beta\in\Delta_i\backslash\{\alpha_i\}$; in particular $G_{\alpha_i}=G_{\alpha_i\beta}$.
  If $g\in G_{\alpha_i}$, then $g$ fixes all of $\Delta_i$; however $G$ acts faithfully on $\Delta_i$ by Proposition~\ref{qtfaith}, thus $g=1$ and $G_{\alpha_i}$ is trivial.
  By definition $t$ must divide $|G_{\alpha_i}|$ and so $t=1$ as well, contradicting that $t>1$.
 \end{proof}
\end{cor}

We now restrict a quasi-transitive group to Theorem~\ref{LPS1}.

\begin{prop}\label{reduct1}
 Let $G$ be a finite quasi-transitive group on a finite set $\Omega$ with orbits $\Delta_1,\dots,\Delta_r$ and $r>1$; then $G$ acts as a finite primitive $\frac{3}{2}$-transitive group on each $\Delta_i$.
 \begin{proof}
  As $d_i=|G_{\alpha_i}|/t$, it is clear that $d_i=|\beta^{G_{\alpha_i}}|$ for any $\beta\in\Delta_i$; hence $G$ is acting $\frac{3}{2}$-transitive on each $\Delta_i$.
  By a classical result about $\frac{3}{2}$-transitive groups, $G$ must act as a Frobenius or primitive group on each $\Delta_i$ \cite[Theorem 10.4]{Wielandt}.
  As $G$ acts faithfully on $\Delta_i$, if $G$ also acts Frobeniusly then $t=|G_{\alpha_i\beta}|=1$ for any $\beta\in\Delta_i$.
  Thus $G$ must act primitively on each $\Delta_i$.
 \end{proof}
\end{prop}

\begin{cor}\label{reduct2}
 Let $G$ be a finite quasi-transitive group on a finite set $\Omega$ with orbits $\Delta_1,\dots,\Delta_r$ and $r>1$; then $G$ acts as a almost simple $\frac{3}{2}$-transitive group on each $\Delta_i$.
 \begin{proof}
  $G$ acts as a primitive $\frac{3}{2}$-transitive group on each $\Delta_i$ by Proposition~\ref{reduct1}.
  Hence $G$ must be either an affine group or an almost simple group by Theorem~\ref{LPS1} .
  Affine groups are characterised by having a unique minimal normal subgroup which is elementary abelian, in which case $G$ is transitive on $\Omega$ by Theorem~\ref{Camina2}.
  Thus $G$ must be almost simple.
 \end{proof}
\end{cor}

\section{The main theorem}

In the previous section we established that if a finite group $G$ acts quasi-transitively on a finite subset $\Omega$ with orbits $\Delta_1,\dots,\Delta_r$ and $r>1$; then $G$ must act as an almost simple $\frac{3}{2}$-transitive group on each $\Delta_i$.
In particular, it must be one of the groups given in Theorem~\ref{LPS2}:
\begin{enumerate}
  \item $G$ is $2$-transitive on $\Omega$.
  \item $n=21$ and $G$ is $A_7$ or $S_7$ acting on the set of pairs of elements of $\{1,\dots,7\}$; the size of the nontrivial subdegrees is 10.
  \item $n=\frac{1}{2}q(q-1)$, where $q=2^f\geq 8$, and either $G=PSL_2(q)$ or $G=P\Gamma L_2(q)$ with $f$ a prime; the size of the nontrivial subdegrees is $q+1$ or $f(q+1)$, respectively.
 \end{enumerate}
As $r>1$, the group $G$ must act on each orbit as one of the above group actions.
Choose two orbits, say $\Delta_1$ and $\Delta_2$; it is enough to show that $G$ can not act on both of these orbits via any combination of the actions given above.

\paragraph{\underline{\bf{Step 1: $G$ does not act as a $2$-transitive group on both $\Delta_1$ and $\Delta_2$}}}

If $G$ acts $2$-transitively on $\Delta_i$, then by applying the Orbit-Stabiliser Theorem first to $G$ and then to $G_{\alpha_i}$, we obtain $|G|=|\Delta_i|(|\Delta_i|-1)\cdot t$.
Hence if $G$ is $2$-transitive on $\Delta_i$, then $|\Delta_i|$ must solve the quadratic equation $tx^2-tx-|G|=0$; however this equation can have at most one positive integer solution.
$G$ can be $2$-transitive on at most one orbit because the $|\Delta_i|$ are pairwise distinct by Corollary~\ref{distinct}.

\paragraph{\underline{\bf{Step 2: $G$ must act as a $2$-transitive group on one of $\Delta_1$ and $\Delta_2$}}}

It is enough to show that $G$ can not act as in ($2$) or ($3$) on the two orbits.
However in both lists, the isomorphism type of $G$ is associated to a unique action.
Hence $G$ can only act on one orbit as a group in ($2$) or ($3$); thus must act on the other orbit as a $2$-transitive group.

In fact we have established that a quasi-transitive group which is not transitive must have exactly $2$ orbits.

\paragraph{\underline{\bf{Step 3: $G$ can not act on $\Delta_i$ as in ($3$)}}}

Assume $G$ acts on $\Delta_1$ as in ($3$); then $G$ is either $PSL_2(q)$ or $P\Gamma L_2(q)$ for $q=2^f\geq 8$.
Also $|\Delta_1|=\frac{1}{2}q(q-1)$ and $d_1=q+1$ or $f(q+1)$ respectively.
As $|G_{\alpha\beta}|=t=\frac{|G|}{|\Delta_1|\cdot d_1}$, it follows that $t=2$ and $G_{\alpha\beta}$ is abelian; thus $G$ is transitive on $\Omega$ by Theorem~\ref{Camina2}.

\paragraph{\underline{\bf{Step 4: The final step}}}

Without loss of generality assume that $G$ acts on $\Delta_1$ as in ($2$) above and acts $2$-transitively on $\Delta_2$.
As $G$ acts as in ($2$), we know $G$ must be $A_7$ or $S_7$ with $|\Delta_1|=21$ and $d_1=10$.
Thus $d_2\cdot|\Delta_2|=d_1\cdot|\Delta_1|=210$ by Corollary~\ref{distinct}.
As $G$ is $2$-transitive on $\Delta_2$, it follows that $d_2=|\Delta_2|-1$ and so $|\Delta_2|\cdot (|\Delta_2|-1)=210$.
In particular, $|\Delta_2|$ must be a positive integer solution to the quadratic equation $x^2-x-210=0$.
This polynomial factors as $x^2-x-210=(x-15)(x+14)$ and so $|\Delta_2|=15$ and $d_2=14$.
Now $d_1=10$ and $d_2=14$ which are not coprime, contradicting Lemma~\ref{dCoprime}.
Hence $G$ can not act $2$ transitively on $\Delta_2$.

\paragraph{}
Thus we have proven the following theorem.
\begin{thm}
 If $G$ is a finite quasi-transitive permutation group on a finite set $\Omega$, then $G$ is transitive on $\Omega$.
\end{thm}

\bibliographystyle{plain}
\bibliography{bibfile}

\end{document}